\documentclass{amsart}
\usepackage{latexsym}
\usepackage{graphics}

\renewcommand{\r}{\rightarrow}
\newcommand{\e}{\epsilon}

\renewcommand{\d}{\delta}
\newcommand{\G}{\Gamma}

\newcommand{\uu}{\underline u}
\newcommand{\ou}{\overline u}
\newcommand{\us}{\underline s}
\newcommand{\os}{\overline s}

\begin{document}

\title[Periodic microstructure]{Asymptotically periodic
$L^2$ minimizers \\ in strongly segregating diblock copolymers}

\author{Adam Chmaj}
\address{Division of Integral Equations \\
Faculty of Mathematics and Information Science \\
Warsaw University of Technology\\ Pl. Politechniki 1 \\
00-661 Warsaw \\ Poland} \email{A.Chmaj@mini.pw.edu.pl}
\thanks{$^*$ A preliminary version of this result was announced at the
``Defects and their dynamics'' workshop (09.08-16.08.2003), Banff
International Research Station, Alberta, Canada.}

\begin{abstract}
Using the delta correction to the standard free energy \cite{bc} in
the elastic setting with a quadratic foundation term and some
parameters, we introduce a one dimension only model for strong
segregation in diblock copolymers, whose sharp interface periodic
microstructure is consistent with experiment in low temperatures.
The Green's function pattern forming nonlocality is the same as in
the Ohta-Kawasaki model. Thus we complete the statement in \cite[p.
349]{rt}$^*$: ``The detailed analysis of this model will be given
elsewhere. Our preliminary results indicate that the new model
exhibits periodic minimizers with sharp interfaces.'' We stress that
the result is unexpected, as the functional is not well posed,
moreover the instabilities in $L^2$ typically occur only along
continuous nondifferentiable ``hairs''.

We also improve the derivation done by van der Waals and use it and
the above to show the existence of a phase transition with Maxwell's
equal area rule. However, this model does not predict the universal
critical surface tension exponent, conjectured to be $11/9$.
Actually, the range $(1.2,1.36)$ has been reported in experiments
\cite[p. 360]{j2}. By simply taking a constant kernel, this exponent
is $2$. This is the experimentally ($ \pm 0.1$) verified tricritical
exponent, found e.g., at the consolute $0.9$ K point in mixtures of
${}^3$He and ${}^4$He. Thus there is a third unseen phase at the
phase transition point.

\end{abstract}

\maketitle

Key to some recent developments in modern microelectronics has been
the ability to create high-quality atomically abrupt interfaces
between different semiconductors, which can produce new quantum
states and uncover unexpected phenomena (see
http://en.wikipedia.org/wiki/Nanotechnology). Several experiments
contradicting common knowledge about bulk materials having only
diffuse interfaces have been reported, see the References,
especially \cite{b}, where a new experimental technique is
explained.

It is remarkable that chemically diverse diblock copolymers and
microemulsions admit very similar characteristic structures, e.g.,
the bicontinuous ordered double diamond \cite{om}. A theory formally
predicting the phase transition between weak and strong segregation
was derived in \cite{bc}. Here we show how to use elasticity and
choose parameters in the Vaserstein pseudoassociation potential to
get local $L^2$ minimizers resembling the square wave \cite{trm}. As
a by-product we get a rigorous proof of the gas to liquid phase
transition.

To make this note self-consistent, we review the history of this
modeling.

The ideal gas law
\begin{equation}
P=\frac{RT}{V}, \label{es}
\end{equation}
where $P$ is the internal pressure, $V$ the volume per mole, $T$ the
absolute temperature, $R$ the gas constant, does not predict a {\em
phase transition between gas and liquid, defined according to
experiment as two densities coexisting at the same pressure} (see
e.g., the isotherms for carbon dioxide in \cite{s}).

Van der Waals formulated the modified equation of state as
\begin{equation}
P=\frac{RT}{V-b} - \frac{a}{V^2} ,\label{mes}
\end{equation}
where $a$ is the attraction parameter which arises from polarization
of molecules into dipoles and $b$ is the volume enclosed within a
particle (the repulsion parameter). Below the critical temperature
isotherms for (\ref{mes}) have a wiggle \cite{s}, which is
unphysical. Maxwell's construction in which the fluid is taken
around a reversible cycle of states from a logical point of view is
worthless \cite{s}, as states on the wiggle have no meaning.
However, such soft reasoning can give other interesting results
\cite{k,s}.

For uniform systems, the Helmholtz free energy is defined as
$\Psi=U-TS$, where $U$ is the internal energy and $S$ the entropy.
The first law of thermodynamics states $dU=\delta Q -PdV$, where
$\delta Q$ is heat change, the second law for reversible systems
states $\delta Q=TdS$. Thus $d\Psi=dU-TdS-SdT=\delta Q-PdV-TdS-SdT$,
so $P=-\frac{\partial \Psi}{\partial V}$ and $\Psi=-\int PdV$. Since
the isothermal compressibility $\beta$ satisfies $\beta\equiv
-\frac{1}{V}\frac{\partial V}{\partial P} \geq 0$, we get
$\frac{\partial^2 \Psi}{\partial V^2}\geq 0$. From (\ref{mes}) we
get
\begin{equation}
\Psi^{\rm unif}=-\int PdV=-RT\ln (V-b)-\frac{a}{V}, \label{Psi}
\end{equation}
and the convex envelope of (\ref{Psi}) gives Maxwell's equal area
rule for the isotherms of (\ref{mes}). However, this and related
statistical mechanics arguments lack a proof of phase transition as
defined above. To be more precise, suppose that Maxwell's envelope
is on the line segment $[V_1 ,V_2 ]$. Why does the system take on
only $V_1$ and $V_2$, but not other values inside the segment? In
other words, the implicit meaning of arguing this way is that {\em
assuming there is a phase transition, it is determined by the equal
area rule}.

As is often the case in mathematics, we can get a proof of phase
transition with Maxwell's equal area rule by studying a related
larger structure. This was done in \cite[p. 40]{j1}. In a related
work, the author also obtained $4/3$ as the critical exponent
\cite[p. 360]{j2}. Other methods have reported $3/2$ (van der
Waals), $1.38$ \cite{pr} and $(1.21,1.32)$ \cite{y}.

We now improve the derivation of van der Waals.

Let $\rho=\frac{1}{V}$ denote the density. The free energy of a
nonuniform state $\rho(x)$ of a nonreversible process is a
functional of the Helmholtz form $I(\rho)=U(\rho)-TS(\rho)$. The
cumulative free energy $\int_0^1 \rho \Psi^{\rm unif} (\rho)dx$ is:
\begin{equation} \Psi^{\rm nonunif}
(\rho)=\int_0^1 \rho \Psi^{\rm unif} (\rho) = \int_0^1 [-RT\rho \ln
(\rho^{-1}-b)-a\rho^2 ].  \label{tent}
\end{equation}
$U(\rho)$ cannot be only $\int_0^1 (C-a\rho )\rho$, but rather
$-\int_0^1 \int_0^1 J(x-y)\rho (x)\rho(y) dxdy \simeq \int_0^1 \rho
(C -a\rho -\frac{c_2}{2} \rho'')dx$. Since the term $\int_0^1 C\rho$
can be neglected due to the mass constraint, one can conclude that
\begin{equation}
I(\rho)=\Psi^{\rm nonunif} (\rho)=\int_0^1 -\frac{c_2}{2} \rho''
\rho +\int_0^1 [-RT\rho \ln (\rho^{-1}-b)-a\rho^2 ]   \label{vdw}
\end{equation}
Note that $-a \rho^2 -RT\rho(\ln(\rho^{-1} -b))$ is concave
(double-well) for $\rho$ and $T$ satisfying $RT\geq (<) ~~2a\rho
(1-\rho b)^2$. However, note that the nonlocal and the $-\int_0^1
a\rho^2$ terms appeared quite separately, therefore it may not be
entirely convincing that one is an extension of the other. Also, if
$-\int_0^1 a\rho^2$ is replaced by $-\int_0^1 \int_0^1 J(x-y)\rho
(x)\rho(y) dxdy$, the solutions of the EL equation are always
continuous, so we have no phase transition as defined above.

In his Nobel lecture, van der Waals expressed his absolute
conviction that molecules associate in complexes not of chemical
origin. He called them pseudoassociations. Thus the derivation can
be improved by adding the nonlocal term to (\ref{tent}):
\begin{equation}
I^{\rm mod}(\rho) =-\int_0^1 \int_0^1 J(x-y)\rho (x)\rho(y) dxdy
+\int_0^1 [-RT\rho \ln (\rho^{-1}-b)-a\rho^2 ] . \label{corr}
\end{equation}
The mass constraint makes possible an addition of a linear term,
which we choose so that $-RT\rho \ln (\rho^{-1}-b)-a\rho^2$ has
equal depth wells. Also, on a bounded interval boundary effects may
come into play, so $J$ is not necessarily translationally invariant:
$J=J(x,y)$. Now (\ref{corr}) is qualitatively the same as derived in
\cite{bc}
\begin{equation}
I(u)=\frac{1}{4} \int_0^1 \int_0^1 J(x,y)(u(x)-u(y))^2 dxdy
+\int_0^1 W(u(x))dx, \label{fe}
\end{equation}
where $W$ is the double-well function
\begin{equation}
W(u)=-\frac{1}{2}ju^2 -\frac{1}{2} u^2 + kT[(1+u)\ln (1+u)+(1-u)\ln
(1-u)], \label{W}
\end{equation}
with $j(x)=\int_0^1 J(x,y)dy$ and $k$ the Boltzmann's constant. For
some works on this and related models see the References, especially
the collective diffusion kinetics of the phase transition in polymer
gels, what the authors in \cite{mt} consider one of the most
exciting problems in current condensed matter physics, and a
nonlocal in time evolution discussed in \cite{tv}.

Let $G(u)=-\frac{1}{2} u^2 + kT[(1+u)\ln (1+u)+(1-u)\ln (1-u)]$,
$G^*$ be the convex envelope of $G$ and
\begin{equation}
I^* (u)=-\frac{1}{2} \int_0^1 J[u](x)u(x)dx +\int_0^1 G^* (u(x)) dx.
\label{conv}
\end{equation}

Let $g^* = {G^*}' $, $[\uu ,\ou ]$ be the interval on which $g^*$ is
constant and $v^* =g^* (\uu )=g^* (\ou )$. We define
\[ \us (v)=\left\{ \begin{array}{ll} {g^*}^{-1} (v), & v \neq  v^*
,\\ \uu , & v=v^* , \end{array} \right. ~~ \os  (v)=\left\{
\begin{array}{ll} {g^*}^{-1} (v), & v \neq  v^* ,\\ \ou , & v=v^* .
\end{array} \right. \]
Let $u$ be a local minimizer of (\ref{fe}). Let $v=g(u)$ and $x_0$
be such that $v(x_0 )=v^*$. We now show that under some conditions
$v'(x_0 )\neq 0$, which implies that the set of discontinuities of
$u$ is finite.

Let $\d_0 >0$ be such that $J(x,y)>\frac{1}{2}J(x,x)$ for all
$(x,y)$ such that $|x-y|<2\d_0$. Let $\d \in (0,\d_0 )$ and $I_\d
=[x_0 -\d ,x_0 +\d ]$. We define
\[ v^+ =\d^2 +{\rm max} \{ v(x)|x\in I_\d \},~~ v^- ={\rm min}\{
v(x)|x\in I_\d \} -\d^2 .\] Let

\begin{equation} \phi = (u-\us (v^-
))\chi_{I_\d} -\frac{1}{1-2\d} \bigl( \int_{I_\d} (u-\us (v^- )
\bigr) \chi_{[0,1]\setminus I_\d}.\label{phi}
\end{equation}

Since $\int_0^1 \phi=0$, we have $I(u)\leq I(u-\phi )$, or
\[ 0\geq \int_0^1 \bigl\{ G^* (u)-G^* (u-\phi )-\phi J[u]+\frac{1}{2}\phi
J[\phi] \bigr\} . \]
On $I_\d$ we have
\[ G^* (u)-G^* (u-\phi )-v\phi \geq (g^* (u-\phi )-v)\phi \geq (v^-
- v^+ )\phi \] and
\[ \frac{1}{2} \int_{I_\d} \phi J[\phi ] \geq \frac{J(x_0 ,x_0 )}{4}
\bigl( \int_{I_\d} \phi \bigr)^2 -\frac{\max_{x\in I_\d}
j(x)}{2(1-2\d)} \bigl( \int_{I_\d} \phi \bigr)^2 \] On
$[0,1]\setminus I_\d$ we have
\[ G^* (u)-G^* (u-\phi )-v\phi \geq (g^* (u-\phi )- g^* (u))
\frac{-1}{1-2\d} \int_{I_\d} \phi \]
\[ \geq \frac{-g'(u)}{(1-2\d )^2} \bigl( \int_{I_\d}
\phi \bigr)^2 + O\Bigl( \bigl( \int_{I_\d} \phi \bigr)^3 \Bigr)
\] and
\[ \frac{1}{2} \int_{[0,1]\setminus I_\d} \phi J[\phi ] \geq
-\frac{\max_{x\in I_\d} j(x)}{2(1-2\d)}  \bigl( \int_{I_\d} \phi
\bigr)^2 + \frac{\int_{[0,1]\setminus I_\d} \int_{[0,1]\setminus
I_\d} J}{2(1-2\d)^2}  \bigl( \int_{I_\d} \phi \bigr)^2
\]
Taking into account these inequalities we obtain
\[ 0 \geq \Bigl\{ \frac{J(x_0 ,x_0 )}{4} + \frac{\int_{[0,1]\setminus I_\d}
\int_{[0,1]\setminus I_\d} J}{2(1-2\d)^2} -\frac{\max_{x\in I_\d}
j(x)}{1-2\d} - \frac{\int_{[0,1]\setminus I_\d} g'(u)}{(1-2\d
)^2}\Bigr\} \bigl( \int_{I_\d} \phi \bigr)^2
\]
\[ + (v^- - v^+ ) \int_{I_\d} \phi + O\Bigl( \bigl( \int_{I_\d} \phi \bigr)^3
\Bigr) . \] Denote the expression in the curly brackets by $C(J,g,\d
)$. After dividing by $\int_{I_\d} \phi$ we get
\[ v^+ -v^- \geq C(J,g,\d ) \int_{I_\d} (u- \us (v^- )) +O(\d^2 ).
\]
In a similar manner, we take
\[ \phi = (u-\os (v^+
))\chi_{I_\d} -\frac{1}{1-2\d} \bigl( \int_{I_\d} (u-\os (v^+ )
\bigr) \chi_{[0,1]\setminus I_\d} \] and obtain
\[ v^+ -v^- \geq C(J,g,\d ) \int_{I_\d} (\os (v^+ )-u) +O(\d^2 ).
\]
Adding these two inequalities and dividing by $2\d$ we get
\[ \frac{v^+ -v^-}{2\d} \geq \frac{C(J,g,\d )}{2} (\os (v^+ )-\us (v^-
)) +O(\d ). \] Letting $\d \r 0$ we get
\[ |v'(x_0 )| \geq \frac{C(J,g)}{2} (\ou -\uu ),\]
where
\[ C(J,g)=\frac{J(x_0 ,x_0 )}{4} + \frac{\int_0^1
\int_0^1 J}{2} - j(x_0 ) - \int_0^1 g'(u) .\] Thus $v'(x_0 )\neq 0$
if $C(J,g)>0$. For now this shows that any local minimizer does not
take values in $(\uu ,\ou )$ and proves the existence of phase
transition.

The well-balanced scaling
\begin{equation}
J_\e (x,y)=\frac{1}{\e} J^s \Bigl( \frac{x-y}{\e} \Bigr) -\e J^l
(x,y), \label{J}
\end{equation}
proposed in \cite{rt}, where $J^s \geq 0$ and $\int_R |x|J^s (x)dx
<\infty$, can be justified as representing an interplay of
attractive ($J^s$) and repulsive ($J^l$) chemical forces. This is
consistent with an established view in physical chemistry, e.g.,
``But the establishment of a well-defined periodicity between the
lamellae, or the rods, depends on the existence of long-range forces
between them: attractive (van der Waals) or repulsive
(electrostatic, steric, plus the short-range Marcelja repulsion)''
\cite[p.2296]{gt}.

Its qualitative shape of that of the ``mexican hat'', a notion that
comes from mathematical biology, where such kernels are also often
used. However, it is interesting to observe that since we assume
only $\int_R |x|J^s (x)dx<\infty$, $J_\e$ can be nonnegative. It is
then its shape, not sign, that gives pattern formation.

The construction of periodic minimizers with discontinuous
interfaces now follows from a sequence of lemmas. Let $I_\e$ denote
$I$ with $J_\e$ given by (\ref{J}). Then $\e^{-1} I_\e$
$\G$-converges as $\e \r 0$ to a singular limit defined on $L^2
(0,1)$ by
\[
I_0 (u)\equiv \left\{ \begin{array}{ll} c_0 \frac{||Du|| (0,1)}{2} -
\frac{1}{4} \int_0^1 \int_0^1 J^l (x,y) (u(x)-u(y))^2 dxdy & u \in
BV((0,1),\{ -1,1\} ), \\
\infty & {\rm otherwise}, \end{array} \right.
\]
where $c_0 \equiv \inf_{\{u: u(\pm \infty)=\pm 1\}} \int_R \int_R
J^s (x-y)(u(x)-u(y))^2dxdy+ \int_R W(u(x))dx$ and $||Du||$ is the
variation measure of $u$ \cite{ab}. Here $\G$-convergence is defined
as
\[
\begin{array}{l}
1.~ For~ every ~\{ u_\e \} \subset L^2 (0,1)~ with~ \lim_{\e \r 0}
u_\e =u,~ \liminf_{\e \r 0} \e^{-1} I_\e (u_\e )\geq I_0 (u) ;\\
2.~For~ every ~ u\in L^2 (0,1)\cap BV ((0,1),\{ -1,1\} ), ~there~
exists~ a~ family~ \\ \hspace{4mm} \{ u_\e \} \subset L^2 (0,1)~
such~that~ \lim_{\e \r 0} u_\e =u,~and~\limsup_{\e \r 0} \e^{-1}
I_\e (u_\e )\leq I_0 (u).    \end{array}
\]
These two inequalities often go together the compactness property
\[
\begin{array}{l}
3. ~~Let ~\e_n ~be ~a ~sequence~
of~positive~numbers~converging~to~0,~ and~ \{ u_n \} ~a
\\ \hspace{4mm} sequence~in~~L^2 (0,1).~If ~\e_n^{-1} I_{e_n}(u_n )
~is~bounded~above~in~n,~then~\{ u_n \} \\ \hspace{4mm}
is~~relatively~~compact~~in~~L^2(0,1)~~ and~~its~~cluster~~
points~~belong~~to\\ \hspace{4mm}BV((0,1),\{ -1,1\} ).
\end{array} \label{c}
\]
A further related property is that a strict local minimum $u_0$ of
$I_0$ perturbs to a local minimum $u_\e$ of $\e^{-1} I_\e$. This was
shown and used in \cite{ks} to obtain local minimizers on dumbbell
domains. In the proof, by a standard argument, a minimizer $u_\e$ of
$\e^{-1} I_\e$ is first constructed in a small closed ball around
$u_0$. Then using properties $1-3$ it is shown that for small enough
$\e >0$, $u_\e$ lies in the interior of the ball, thus is a local
minimizer of $\e^{-1} I_\e$ in $L^2$.

When $G$ is nonconvex this argument goes through a convexification
and the calculation above. It is shown that a local minimum of
$I^*_0$ in a subspace of $BV ((0,1),\{ -1,1 \} )$ having a fixed
number of jumps is also a local minimum of $I^*_0$ in $L^2 (0,1)$.
This observation reduces the problem of determining local minima of
$I^*_0$ to a finite dimensional one. Since
$\frac{1}{2}(||Du||(0,1))$ is equal to the number of jumps of $u$,
it is enough to investigate only the second term of $I^*_0$. In
general, its critical points are determined from a system of
nonlinear algebraic equations and this is not an easy, or indeed,
pleasant task.

The Ohta-Kawasaki functional follows from \cite{am}, though it has
been shown in \cite{om} that it is not sophisticated enough to
predict the complicated phase separation morphologies. However, to
confirm that pseudoassociations exist and to calculate the critical
point exponent, we make a study in one dimension.

The robust functional
\begin{equation}
I(w)= \int_0^1 [\frac{1}{2} \e^2 w''(x)^2 +W(w'(x)) +w^2 (x) ]dx,
\label{ef}
\end{equation}
can model e.g., the twinned martensite phase in Nitinol \cite{mu}.
The elastic foundation third term in (\ref{ef}) is nonlocal. Namely,
let us consider (\ref{ef}) with the boundary conditions
$w(0)=w(1)=0$. Let $w=v'$, $v=(-D^2)^{-1} (u-m)$, where
\[ -D^2 : \{ v\in W^{2,2} :v'(0)=v'(1)=0, \int_0^1 v=0 \} \r \{ w\in
L^2 : \int_0^1 w=0 \} . \] Then $w'=v''=m-u$, $w''=-u'$ and
\[ \int_0^1 w^2 =\int_0^1 v'^2 =-\int_0^1 v''v =\int_0^1 (u-m)(-D^2
)^{-1} (u-m) =\int_0^1 u (-D^2 )^{-1} u \] and (\ref{ef}) becomes
\[ I(u)=\int_0^1 [\frac{1}{2} \e^2 u'(x)^2 +W(m-u(x))]dx +\int_0^1
\int_0^1 G(x,y) u(x)u(y)dxdy ,\] with $G$ as in \cite{rw}. Using the
improved van der Waals derivation, we obtain (\ref{fe}) with
$J=J_\e$ - the strong separation diblock copolymer functional. With
$G$ and small $\e$ - asymptotically periodic local minimizers with
sharp interfaces.

\end{document}